\documentclass[a4paper,10pt,reqno]{amsart}
\usepackage[a4paper,tmargin=40mm,bmargin=35mm,hmargin=30mm]{geometry}
\usepackage[T1]{fontenc}
\usepackage{lmodern}
\usepackage{amsmath}
\usepackage{amssymb}
\usepackage{amsthm}
\usepackage{stmaryrd}
\usepackage{tikz}
\usepackage[all]{xy}
\usepackage{booktabs}
\usepackage{multirow}
\usepackage{multicol}
\usepackage{here}
\usepackage{aliascnt}
\usepackage{hyperref}
\usepackage[noabbrev]{cleveref}

\newcommand{\NewTheorem}[2]{
	\newaliascnt{#1}{TheoremEnvironment}
	\newtheorem{#1}[#1]{#1}
	\aliascntresetthe{#1}
	\crefname{#1}{#1}{#2}
	\Crefname{#1}{#1}{#2}
}

\theoremstyle{definition}

\NewTheorem{Definition}{Definitions}
\NewTheorem{Remark}{Remarks}
\NewTheorem{Axiom}{Axioms}
\NewTheorem{Example}{Examples}
\NewTheorem{Observation}{Observations}
\NewTheorem{Convention}{Conventions}
\NewTheorem{Notation}{Notations}
\NewTheorem{Setting}{Settings}
\NewTheorem{Question}{Questions}
\NewTheorem{Answer}{Answers}
\NewTheorem{Conjecture}{Conjectures}
\NewTheorem{Problem}{Problems}
\NewTheorem{Solution}{Solutions}
\NewTheorem{Goal}{Goals}
\NewTheorem{Comment}{Comments}
\NewTheorem{Aim}{Aims}
\NewTheorem{Caution}{Cautions}
\NewTheorem{Exercise}{Exercises}
\NewTheorem{Examples}{Examples}
\theoremstyle{plain}
\NewTheorem{Proposition}{Propositions}
\NewTheorem{Lemma}{Lemmas}
\NewTheorem{Theorem}{Theorems}
\NewTheorem{Corollary}{Corollaries}

\crefname{enumi}{}{}
\Crefname{enumi}{}{}
\creflabelformat{enumi}{(#2#1#3)}
\crefname{enumii}{}{}
\Crefname{enumii}{}{}
\creflabelformat{enumii}{(#2#1#3)}
\crefname{enumiii}{}{}
\Crefname{enumiii}{}{}
\creflabelformat{enumiii}{(#2#1#3)}

\makeatletter
\renewcommand{\p@enumii}{}
\renewcommand{\p@enumiii}{}
\makeatother

\numberwithin{equation}{section}
\crefname{equation}{}{}
\Crefname{equation}{}{}
\creflabelformat{equation}{(#2#1#3)}

\newcommand{\SwapSymbols}[1]{
	\expandafter\let\expandafter\temporarysymbol\csname #1\endcsname
	\expandafter\let\csname #1\expandafter\endcsname\csname var#1\endcsname
	\expandafter\let\csname var#1\endcsname\temporarysymbol
}

\SwapSymbols{epsilon}
\SwapSymbols{phi}
\SwapSymbols{Gamma}
\SwapSymbols{Delta}
\SwapSymbols{Theta}
\SwapSymbols{Lambda}
\SwapSymbols{Xi}
\SwapSymbols{Pi}
\SwapSymbols{Sigma}
\SwapSymbols{Upsilon}
\SwapSymbols{Phi}
\SwapSymbols{Psi}
\SwapSymbols{Omega}

\newcommand{\bbZ}{\mathbb{Z}}

\newcommand{\cA}{\mathcal{A}}

\newcommand{\cC}{\mathcal{C}}

\newcommand{\cL}{\mathcal{L}}

\newcommand{\cN}{\mathcal{N}}
\newcommand{\cO}{\mathcal{O}}

\newcommand{\cS}{\mathcal{S}}
\newcommand{\cT}{\mathcal{T}}

\newcommand{\cW}{\mathcal{W}}
\newcommand{\cX}{\mathcal{X}}
\newcommand{\cY}{\mathcal{Y}}

\newcommand{\To}{\longrightarrow}

\DeclareMathOperator{\Hom}{Hom}
\DeclareMathOperator{\End}{End}

\DeclareMathOperator{\PSupp}{PSupp}

\DeclareMathOperator{\lAnn}{lAnn}
\DeclareMathOperator{\Mod}{Mod}
\let\mod\relax
\DeclareMathOperator{\mod}{mod}

\DeclareMathOperator{\nPSpec}{nPSpec}
\DeclareMathOperator{\nSupp}{nSupp}
\DeclareMathOperator{\noeth}{noeth}

\DeclareMathOperator{\Zg}{Zg}

\DeclareMathOperator{\nPSupp}{nPSupp}
\DeclareMathOperator{\PSpec}{PSpec}

\DeclareMathOperator{\nsSupp}{nsSupp}

\DeclareMathOperator{\rAnn}{rAnn}

\DeclareMathOperator{\Ann}{Ann}

\let\Im\relax
\DeclareMathOperator{\Im}{Im}
\DeclareMathOperator{\sPsupp}{sSupp}
\DeclareMathOperator{\Spec}{Spec}

\DeclareMathOperator{\Ass}{Ass}
\DeclareMathOperator{\Supp}{Supp}

\DeclareMathOperator{\ASpec}{ASpec}

\DeclareMathOperator{\ASupp}{ASupp}
\DeclareMathOperator{\sPSupp}{sPSupp}
\DeclareMathOperator{\nsPSupp}{nsPsupp}

\title{Spectrum of an abelian category via premonoform objects}
\author{Reza Sazeedeh}

\address{Department of Mathematics, Urmia University, P.O.Box: 165, Urmia, Iran}
\email{rsazeedeh@ipm.ir and r.sazeedeh@urmia.ac.ir}

\begin{document}

\begin{abstract}
Let $\cA$ be an abelian category. In this paper, we study ${\rm (n)PSpec}\cA$, a topological space formed by equivalence classes derived from an equivalence relation on (noetherian) premonoform objects. We classify torsion classes of $\cA$ via closed subclasses of $\nPSpec\cA$.  We introduce a new topology on $\PSpec\cA$ and we classify Serre subcategories of $\noeth A$ and localizing subcategories of $\cA$ using this topology. If $A$ is a commutative noetherian ring, we show that $\nPSpec A$ is homeomorphic to $\Spec A$. Moreover, there is a one-to-one correspondence between the closed subsets of $\nPSpec A$ and the open subsets of $\ASpec A$, the atom spectrum of $A$. Finally, we explore the relationships between the new subctegories of $\Mod A$ and subsets of $\nPSpec A$ introduced in this paper, and the known subcategories of $\Mod A$ and subsets of other spectra of $A$.   
\end{abstract}

\maketitle
\tableofcontents

\section{Introduction}

The spectrum of an abelian category plays a central role in identifying its subcategories. There are several types of spectra depending on the specific category in question. A classical and well-known spectrum is Spec $A$ in the category of $A$-modules, where $A$ is a commutative noetherian ring. Let $\Mod A$ be the category of $A$-modules and mod-$A$ be the subcategory of finitely generated $A$-modules. Gabriel [Ga] established a one-to-one correspondence between localizing subcategories of $\Mod A$, Serre subcategories of $\mod A$ and specialization-closed subsets  of Spec $A$. When $A$ is a commutative ring, the classification of subcategories of $\Mod A$ in terms of subsets of $\Spec A$ has been studied by several authors (e.g. [Hop], [N], [GP], [Ho]).  In the case where $A$ is a commutative noetherian ring, Takahashi [T] established isomorphisms between the lattices of subcategories of $\Mod A$ and subsets of $\Spec A$ posed by Hovey [H].    

 For a locally coherent Grothendieck category $\cA$, the Ziegler spectrum Zg$\cA$ is a topological space formed by the isomorphism classes of indecomposable injective objects (see for example [H], [Kr]). This concept was originally introduced by Ziegler [Z].  Herzog [H] established a one-to-one correspondence between open subsets of Zg$\cA$ and Serre subcategories of the abelian category coh-$\cA$  consisting of coherent objects in $\cA$. 
 
 In an arbitrary abelian category $\cA$, Kanda [K] introduced the atom spectrum $\ASpec \cA$  which is a topological space consisting of atoms. Atoms in an abelian category are defined as the equivalence classes of an equivalence relation on monoform objects, drawing inspiration form monoform modules over noncommutative rings, as explored by Storrer [St]. In the context of a locally noetherian Grothendieck category $\cA$, the topologies Zg$\cA$ and ASpc$\cA$ are homeomorphic.
 
To classify nullity classes in an abelian category $\cA$, Liu and Stanley [LS] introduced  $\PSpec\cA$, a topological space formed by equivalence classes derived from an equivalence relation on premonoform objects, based on their nullity classes. The topology is constructed by defining closed subclasses. The primary aim of this paper is to establish connection between subcategories of $\cA$ and closed subclasses of $\PSpec\cA$. We demonstrate that the subspace $\nPSpec \cA$ consisting of equivalence classes in $\PSpec \cA$ represented by noetherian premonoform objects, plays a central role in classifying the subcategories of $\cA$. As an application of our results, we explore new subcategories of $\Mod A$ and subsets of $\nPSpec A$ in the case where $A$ is a commutative noetherian ring. Furthermore, we extend the lattices isomorphisms established by  Takahashi [T] and the classifications investigated by Stanley and Wang [SW].   

 In Section 2, we assume that $\cA$ is an abelian category and study some basic properties of $\PSpec \cA$. We show that $\PSpec \cA$ is an Alexandroff and Kolmogorove space. One of the main results in this section is the following theorem.  
 
\begin{Theorem}\label{tors}
There exists an injective map $\Phi\mapsto \nPSupp^{-1}\Phi$  between closed and extension-closed subclasses of $\nPSpec\cA$
and torsion classes of $\cA$. This map has a left inverse given by  $\cX\mapsto\nPSupp\cX$.
\end{Theorem}

We also classify  torsion classes of weakly finite type of a locally noetherian Grothendieck category as follows. 
\medskip
\begin{Theorem}\label{classloc}
 Let $\cA$ be a locally noetherian Grothendieck category. Then the map $\cX\mapsto\nPSupp\cX$ establishes a one-to-one correspondence between torsion classes of weakly finite type of $\cA$ and closed and extension-closed subsets of $\nPSpec\cA$. The inverse map is given by $\Phi\mapsto \overrightarrow{\nPSupp^{-1}\Phi\cap\noeth A}.$
\end{Theorem}
We show that if the collection of all closed and extension-closed subclasses of $\nPSpec\cA$ forms a topology on $\nPSpec\cA$, then a closed subclass $\Phi$ of $\nPSpec\cA$ is compact if and only if there exists a noetherian object $M$ in $\cA$ such that $\Phi=\nPSupp M$; see \cref{compp}.

In Section 3, we study premonoform $A$-modules over a (non)commutative ring $A$. Suppose that $A$ is left noetherain and $M$ is an $A$-bimodule  and a premonoform right $A$-module. We show that $\lAnn_AM$ is an extremely prime ideal of $A$ and if $M$ is a uniform left $A$-module, then $M$ is a monoform left $A$-module; see \cref{pp}. An analogues result holds when $A$ is right noetherian and $M$ is an $A$-bimodule that is a premonoform left $A$-module. As a conclusion, over a commutative noetherian ring $A$, every uniform and premonoform $A$-module is monoform. Let $\frak a$ be a right ideal of $A$. Then we show that  $A/\frak a$ is a premonoform  right $A$-module if and only if $\frak a$ is  completely prime; \cref{compr}.

In Section 4, for an abelian category $\cA$, we  introduce a new topological space on $\PSpec \cA$. A closed subclass of $\PSpec\cA$ is called Serre closed and it is defined in terms of Serre subcategories of objects in $\cA$.
Over a commutative noetherian ring $A$,  closed subsets of $\nPSpec A$ are precisely  Serre closed subsets and  extension-closed subsets of $\nPSpec A$ are Serre extension-closed subsets, see \cref{comclex}. One of the main results in this section is the following theorem.

\begin{Theorem}
  There is a one-to-one correspondence between Serre subcategories of $\noeth\cA$ and Serre closed  and Serre extension-closed subclasses of $\nPSpec\cA$. 
  \end{Theorem}
The following theorem classifies localizing subcategories of a locally noetherian Grothendieck category $\cA$ via the new topology on $\nPSpec\cA$. 
\begin{Theorem}
Let $\cA$ be a locally noetherian Grothendieck category. Then the map $\cX\mapsto\nsPSupp\cX$ establishes a one-to-one correspondence between localizing subcategories of $\cA$ and Serre closed and Serre extension-closed subsets of $\nPSpec\cA$.  The inverse map is given by $\Phi\mapsto\nsPSupp^{-1}\Phi$.
 \end{Theorem}

In Section 5, we study premonoform $A$-modules over a commutative ring $A$. We show that if $A$ is noetherian, then closed and extension-closed subsets of $\nPSpec\cA$ are in correspondence with specialization-closed subsets of $\Spec A$; see \cref{spec}. An immediate conclusion of this result implies that every torsion class of $\mod A$ is Serre. For every premonoform module $H$, we set $\frak p_H=\Ann_AH$ which is a prime ideal of $A$. One of the main results in this section as follows.

\begin{Theorem}\label{homeom}
The function $p:\PSpec A\To\Spec A$; given by $[H]\mapsto \frak p_H$ is continuous. Moreover, if $A$ is notherian, then $p:\nPSpec A\To\Spec A$ is homeomorphism.  
\end{Theorem}
As a conclusion of the above theorem, every closed subset of $\nPSpec A$ is extension-closed. 
 The other main result in this section is to establish a relation between closed subsets of $\nPSpec A$ and open subsets of $\ASpec A$. We have the following theorem.
  \begin{Theorem}
Let $A$ be a noetherian ring. Then there exists a bijective map $\theta:\nPSpec A\to \ASpec A$ given by $[H]\mapsto \overline{H}$. Moreover, $\theta$ establishes a one-to-one correspondence between closed subsets of $\nPSpec A$ and open subsets of $\ASpec A$.    
\end{Theorem}
  
  As the last conclusion, we provide the following classifications.
  
\begin{Theorem}
Let $A$ be a commutative noetherian ring and consider the following lattices.
 
$\bullet$ $\cL_{\rm thick}(D_{\rm perf}(A))$, the lattice of all thick subcategories of $D_{\rm perf}(A)$;
 
 $\bullet$ $\cL_{\rm null}(\mod A)$, the lattice of all nullity classes of $\mod A$;
 
 $\bullet$ $\cL_{\rm coh}(\mod A)$, the lattice of all coherent subcategories of $\mod A$;
 
 $\bullet$ $\cL_{\rm Serre}(\mod A)$, the lattice of all Serre subcategories of $\mod A$;
 
  $\bullet$ $\cL_{\rm narrow}(\mod A)$, the lattice of all narrow subcategories of $\mod A$;
 
 $\bullet$ $\cL_{\rm loc}(\Mod A)$, the lattice of all localizing subcategories of $\Mod A$;
 
  $\bullet$ $\cL_{\rm torsft}(\Mod A)$, the lattice of all torsion classes of finite type of $\Mod A$;
 
 $\bullet$ $\cL_{\rm torswft}(\Mod A)$, the lattice of all torsion classes of weakly finite type of $\Mod A$;
 
$\bullet$ $\cL_{\rm spcl}(\Spec A)$, the lattice of all specialization-closed subsets of $\Spec A$;

$\bullet$ $\cL_{\rm clos}(\nPSpec A)$, the lattice of all closed subsets of $\nPSpec A$;

$\bullet$ $\cL_{\rm cexc}(\nPSpec A)$, the lattice of all closed and extension-closed subsets of $\nPSpec A$;

$\bullet$ $\cL_{\rm sclos}(\nPSpec A)$, the lattice of all Serre closed subsets of $\nPSpec A$;

$\bullet$ $\cL_{\rm scsexc}(\nPSpec A)$, the lattice of all Serre closed and Serre extension-closed subsets of $\nPSpec A$;

$\bullet$ $\cL_{\rm open}(\ASpec A)$, the lattice of all open subsets of $\ASpec A$;

$\bullet$ $\cL_{\rm Zieg}(\Zg A)$, the lattice of all open subsets of $\Zg A.$

Then there are the following equalities and isomorphisms of lattices:
$$ \cL_{\rm Zieg}(\Zg A)\cong\cL_{\rm open}(\ASpec A)
 \cong \cL_{\rm clos}(\nPSpec A)
  =\cL_{\rm cexc}(\nPSpec A)=\cL_{\rm sclos}(\Spec A)$$$$
  =\cL_{\rm cexc}(\nPSpec A)\cong \cL_{\rm spcl}(\Spec A)\cong \cL_{\rm loc}(\Mod A)=\cL_{\rm torsft}(\Mod A)
  =\cL_{\rm torswft}(\Mod A)$$$$\cong \cL_{\rm Serre}(\mod A)
  =\cL_{\rm null}(\mod A)\\=\cL_{\rm narrow}(\mod A)=\cL_{\rm coh}(\mod A)\cong \cL_{\rm thick}(D_{\rm perf}(A)).$$
 \end{Theorem}

\section{Premonoform objects in an abelian category}
In this section, we assume that $\cA$ is an abelian category. We briefly list some necessary basic definitions and notations.

Let $\cS$ be a subcategory of $\cA$. We denote by $\cS_{\rm sub}$ ($\cS_{\rm quot}$), the smallest subcategory of $\cC$ containing $\cS$ which is closed under taking subobjects (quotients). These subcategories can be specified as follows:    
$$\langle\cS\rangle_{\rm sub}=\{N\in\cA|\hspace{0.1cm} N \hspace{0.1cm}{\rm is\hspace{0.1cm} a\hspace{0.1cm} subobject\hspace{0.1cm} of\hspace{0.1cm} an \hspace{0.1cm}object\hspace{0.1cm} in}\hspace{0.1cm} \cS\};$$
$$\langle\cS\rangle_{\rm quot}=\{M\in\cA|\hspace{0.1cm} M \hspace{0.1cm}{\rm is\hspace{0.1cm} a\hspace{0.1cm} quotient\hspace{0.1cm}
object\hspace{0.1cm} of\hspace{0.1cm} an \hspace{0.1cm}object\hspace{0.1cm} in}\hspace{0.1cm} \cS\}.$$

Let $\cT$ be another subcategory of $\cA$. We denote by $\cS\cT$, the {\it extension subcategory} of $\cS$ and $\cT$ which is:
$$\cS\cT=\{M\in\cA|\hspace{0.1cm} {\rm there \hspace{0.1cm} exists \hspace{0.1cm} an\hspace{0.1cm}exact\hspace{0.1cm}
 sequence}$$$$
0\To L\To M\To N\To 0\hspace{0.1cm} {\rm with}\hspace{0.1cm} L\in\cS
\hspace{0.1cm}{\rm and}\hspace{0.1cm} N\in\cT\}.$$

For any $n\in \mathbb{N}_0$, we set $\cS^0=\{0\}$ and
$\cS^n=\cS\cS^{n-1}$. In the case where $\cS^2=\cS$, we say
that $\cS$ is {\it closed under taking extension}. We also define
$\langle\cS\rangle_{\rm ext}=\bigcup_{n\geq 0}\cS^n$, the smallest
subcategory of $\cA$ containing $\cS$ which is closed under taking
extensions.

\begin{Definition} We recall some the definition of some well known subcategories of $\cA$.

$\bullet$   A full subcategory $\cN$ of $\cA$ is said to be a {\it nullity class} if it is closed under taking quotients and extensions.

$\bullet$ A full subcategory $\cT$ of $\cA$ is said to be a {\it torsion class} if it is a nullity class and it is closed under arbitrary direct sums.
 
 $\bullet$ A full subcategory $\cC$ of $\cA$ is said to be {\it narrow} if it is
closed under taking cokernels and extensions.
 
 $\bullet$ A full abelian subcategory $\cW$ of $\cA$ is said to be {\it coherent (or wide)} if it is closed under taking extensions.
 
$\bullet$ A full subcategory $\cS$ of $\cA$ is said to be {\it Serre} if it is
closed under taking subobjects, quotients and extensions.

$\bullet$ A full subcategory $\cX$ of $\cA$ is said to be {\it localizing} if $\cX$ is a Serre subcategory closed under taking arbitrary direct sums.  
\end{Definition}

 For any subcategory $\cX$ of $\cA$,  $\sqrt{\cX}$ is the smallest Serre subcategory of $\cA$ containing $\cX$. Clearly $\sqrt{\cX}=\langle\langle\langle\cX\rangle_{\rm sub}\rangle_{\rm quot} \rangle_{\rm ext}$. The localizing subcategory generated by $\cX$ is denoted by $\langle\cX\rangle_{\rm loc}$.
 
 \medskip

\begin{Lemma}\label{kan}
Let $\cS$ be a subcategory of $\cA$. Then
\begin{center}
 $\langle\cS\rangle_{\rm ext}=\{N\in\cA|\hspace{0.1cm}\textnormal{there exists a finite filtration}\hspace{0.1cm} 0=N_0\subset N_1\subset N_2\subset\dots \subset N_n=N\hspace{0.1cm}\textnormal{of subobjects of }\hspace{0.1cm} N\hspace{0.1cm}\textnormal{such that}\hspace{0.1cm} N_i/N_{i-1}\in\cS\hspace{0.1cm}\textnormal{for each}\hspace{0.1cm} 1\leq i\leq n\}.$
 \end{center}
\end{Lemma}
\begin{proof}
 Given $N\in\langle\cS\rangle_{\rm ext}$, there exists a positive integer $n$ such that $N\in\cS^n$. We prove by induction on $n$ that there exists a filtration 
 $0\subset N_1\subset N_2\subset\dots \subset N_n=N$ of subobjects of $N$ of length $n$ such that $N_i/N_{i-1}\in\cS$ for each $1\leq i\leq n$. If $n=1$, there is nothing to prove. Since $N\in\cS^n$, there exists an exact sequence $0\To N_1\To N\To N/N_1\To 0$ such that $N_1\in\cS$ and $N/N_1\in\cS^{n-1}$. By the induction hypothesis, there exists a finite filtration  $0\subset N_2/N_1\subset \dots \subset N_n/N_1=N/N_1$ such that $N_i/N_{i-1}\in\cS$ for each $2\leq i\leq n$. Consequently, $0\subset N_1\subset N_2\subset\dots\subset N_n=N$ is the desired filtration.  
\end{proof}

\medskip
\begin{Definition}
A nonzero object $M$ in $\cA$ is {\it monoform} if for every nonzero subobject $N$ of $M$, there is no common nonzero subobjects of $M$ and $M/N$.  We denote by $\ASpec _0\cA$, the set of all monoform
objects of $\cA$. Two monoform objects $H$ and $H'$ in $\ASpec _0\cA$ are said to be {\it
atom-equivalent} if they have a common nonzero subobject. The atom equivalence establishes an
equivalence relation on monoform objects; and hence for every
monoform object $H$, we denote the
 {\it equivalence class} of $H$, by $\overline{H}$, that is
\begin{center}
  $\overline{H}=\{G\in\ASpec _0\cA|\hspace{0.1cm} H \hspace{0.1cm}
 {\rm and}\hspace{0.1cm} G$ has a common nonzero subobject$\}.$
\end{center}
The {\it atom spectrum} $\ASpec\cA$ of $\cA$ is the quotient set
of $\ASpec _0\cA$ consisting of all equivalence classes induced by
this equivalence relation; in other words 
$\ASpec\cA=\{\overline{H}|\hspace{0.1cm} H\in\ASpec _0\cA\}.$
\end{Definition}
\medskip

\begin{Definition}
 A nonzero object $M$ in $\cA$ is {\it premonoform} if it does not have nontrivial quotient $M/N$ as subobjects of $M$. Clearly, every monoform object is premonoform.  
We denote by $\PSpec_0\cA$, the collection of isomorphism classes of premonoform objects in $\cA$. For a subcategory $\cX$ of $\cA$, we denote by $\overline{C(\cX)}$, the nullity class generated by $\cX$ that is the smallest nullity class containing $\cX$. It is clear that $\overline{C(\cX)}=\langle\langle \cX\rangle_{\rm quot}\rangle_{\rm ext}$. If $\cX=\{X\}$ has one element, we set $\overline{C(\cX)}=\overline{C(X)}$.
\end{Definition}

\medskip
\begin{Definition}
Assume that $A,B\in\PSpec_0\cA$. We define  $A\sim B$ if and only if $\overline{C(A)}=\overline{C(B)}$, i.e. they generate the same nullity class.
\end{Definition}
 
 Clearly $\sim$ is an equivalence relation on $\PSpec_0\cA$. The {\it pspectrum} of $\cA$ is denoted by $\PSpec\cA=\PSpec_0\cA/\sim$ which is the collection of the equivalence classes $[H]$ of premonoform objects $H$ in $\PSpec_0\cA$. If $\cA$ is skeletally small, then $\PSpec\cA$ is a set.  For every object $M$ in $\cA$, the {\it psupport of} $M$ is defined as follows 
  $$\PSupp M=\{[H]\in\PSpec\cA|\hspace{0.1cm} H\in\overline{C(M)}\}=\{[H]\in\PSpec\cA|\hspace{0.1cm} [H]\subset\overline{C(M)}\}.$$
  
  The psupport of a subcategory $\cX$ of $\cA$ is defined as $\PSupp \cX=\bigcup_{M\in\cX}\PSupp M$ and for any subclass $\Phi$ of $\PSpec\cA$ we define $\PSupp^{-1}\Phi=\{M\in\cA|\hspace{0.1cm}\PSupp M\subset \Phi\}$.   A subclass $\Phi$ of $\PSpec \cA$ is said to be {\it closed} if for every $[H]\in\Phi$, we have $\PSupp H\subset\Phi$. The subclass $\Phi$ is said to be {\it extension-closed} if $\PSupp^{-1}\Phi$ is closed under extensions. We also have the following lemma:
  
 \medskip
 \begin{Lemma}\label{sursum}
  Let $\cX$ be  a subcategory of $\cA$ closed under finite direct sums. Then
  $$\PSupp\cX=\{[H]\in\PSpec\cA|\hspace{0.1cm} H\in\overline{C(\cX)}\}.$$
   \end{Lemma}
   \begin{proof}
  Given $[H]\in\PSupp\cX$, there exists $M\in\cX$ such that $[H]\in\PSupp M$. Then $H\in\overline{C(M)}\subset\overline{C(\cX)}$. Conversely if $[H]$ is an element in $\PSpec\cA$ such that $H\in\overline{C(\cX)}$, then by \cref{kan}, $H$ has a finite filtration of subobjects 
  $$0\subset H_1\subset H_2\subset\dots H_{n-1}\subset H_n=H$$ and $X_i\in\cX$ such that $H_i/H_{i-1}$ is a quotient of $X_i$ for $i=1,\dots n$. By the assumption, $X=\bigoplus_{i=1}^n X_i$ is in $\cX$ and it is clear that $H_i/H_{i-1}$ is a quotient of $X$ for each $i$. Consequently, $[H]\in\PSupp X$.

     \end{proof}

\begin{Lemma}
Let $M$ be a finitely generated object in $\cA$. Then $\PSupp M$ is a non-empty class.
\end{Lemma}
\begin{proof}
It is clear that $M$ has a maximal subobject $N$ and so $[M/N]\in\PSupp M$. 
\end{proof}

 A topological space $X$ is called {\it Alexandroff} if the union of any family of closed subclass of $X$ is closed. A topological space $X$ is said to be {\it Kolmogorov} (or {\it $T_0$-space}) if for any distinct points $x,y$ of $X$, there exists an closed subclass of $X$ containing exactly one of them.

 \begin{Proposition}\label{kolmo}
The topological space  $\PSpec\cA$ is Alexandroff and Kolmogorov. 
 \end{Proposition}
 \begin{proof}
 It is clear by definition that $\PSpec\cA$ is Alexandroff. Let $[H]$ and $[G]$ be distinct elements of $\PSpec\cA$ such that there is no closed subclass of $\PSpec\cA$ containing exactly one of them. Then $\overline{G}\in\PSupp H$ and so $\overline{C(G)}\subset\overline{C(H)}$. Symmetrically $\overline{C(H)}\subset\overline{C(G)}$ and so $\overline{C(G)}=\overline{C(H)}$. Thus $[H]=[G]$ which is a contradiction.
 \end{proof}
 
 As monoform objects are premonoform, for a monoform object $H$ in $\cA$, the relationship between the atom $\overline{H}\in\ASpec\cA$ and  the equivalence class $[H]\in\PSpec\cA$ is of particular importance. We have the following lemma.
 \medskip
 
 \begin{Lemma}\label{isosub}
 Let $H$ and $G$ be monoform objects in $\cA$ such that $\overline{G}=\overline{H}$ and $[G]=[H]$. Then $G$ is isomorphic to a subobject of $H$.
 \end{Lemma}
 \begin{proof}
 As $[H]\in\PSupp G$, we have $H\in\overline{C(G)}$. Then there exists a monoform object $G_1$ in $\cA$ isomorphic to $ G$ such that $H$ contains some quotient $G_1/K$ of $G_1$. But the condition $\overline{G}=\overline{H}$ implies that $K=0$ and so $G_1$ is a subobject of $H$. 
 \end{proof}
 
\medskip

\begin{Example}
Consider the abelian gropus $\mathbb{Q}$ and $\mathbb{Z}$. It is clear that they are monoform and $\overline{\mathbb{Q}}=\overline{\mathbb{Z}}$. But since $\mathbb{Q}$ is not isomorphic to a abelian subgroup of $\mathbb{Z}$, it follows from \cref{isosub} that 
$[\mathbb{Z}]\neq [\mathbb{Q}]$.
\end{Example}

 \begin{Lemma}
 If $[H]\in\PSpec\cA$ has a noetherian representative, then  $H$ is noetherian.
 \end{Lemma}
 \begin{proof}
 If $G$ is a  noetherian premonoform object in $\cA$ such that $[H]=[G]$, then we have $H\in\overline{C(G)}$. Therefore $H$ is noetherian by \cref{kan}.
 \end{proof}

 \medskip
 We denote by $\nPSpec\cA$, the class of all elements of $\PSpec\cA$ represented by noetherian premonoform objects.  For every object $M$ of $\cA$, we set $\nPSupp M=\nPSpec\cA\cap \PSupp M$. If $M$ is a noetherian object in $\cA$, it is clear that $\nPSupp M=\PSupp M$. For every subcategory $\cX$ of $\cA$, we set $\nPSupp\cX=\bigcup_{M\in\cX}\nPSupp M$. We denote by $\noeth A$ the abelian subcategory of $\cA$ consisting of noetherian objects. It is clear that $\PSpec({\rm noeth}\cA)=\nPSpec\cA$. In particular, $\nPSpec\cA$ is a topological subspace of $\PSpec\cA$. It is well-known set theoretic observation  that if $\cA$ is a locally noetherian category, then $\noeth \cA$ is skeletally small (see for example [K, Proposition 5.2]). Then in this case $\nPSpec \cA$ is a set.

\medskip
\begin{Lemma}\label{exttt}
Let $M$ be a noetherian object in $\cA$. Then $\nPSupp M$ is an extension-closed subclass of $\nPSpec\cA$.
\end{Lemma}
\begin{proof}
Since $M$ is noetherian, all objects in $ \overline{C(M)}$ are noetherian by \cref{kan}.   It suffices to show that $\PSupp M=\PSupp\overline{C(M)}$ and so the assertion is attained by [LS, Theorem 6.8].
Clearly $\PSupp M\subset \PSupp \overline{C(M)}$. Conversely, if $[H]\in\PSupp \overline{C(M)}$, then there exists an object $X$ in $\overline{C(M)}$ such that $[H]\in\PSupp X$. Then $H\in \overline{C(X)}\subset \overline{C(M)}$ so that $[H]\in\PSupp M$.
\end{proof} 
 \medskip

 \begin{Lemma}\label{filt}
Let $M$ be a noetherian object in $\cA$. Then $M$ has a finite filtration of subobjects
$$0=M_0\subset M_1\subset \dots \subset M_{n-1}\subset M_n=M$$ such that $M_i/M_{i-1}$ is a quotient of some noetherian premonoform object  $H_i\in\cA$ with $[H_i]\in\PSupp M$ for $i=1,\dots,n$. 
 \end{Lemma}
 \begin{proof}
 Consider $\cS=\{H\in\cA|\hspace{0.1cm} H \hspace{0.1cm}\textnormal{is a premonoform object in}\hspace{0.1cm}\cA\hspace{0.1cm} \textnormal{such that}\hspace{0.1cm} [H]\in\PSupp M\}$. Since $M$ is noetherian, each object in $\cS$ is noetherian. Clearly $\PSupp M\subset \PSupp\overline{C(\cS)}$. Then $M\in\overline{C(\cS)}$ by [LS, Proposition 6.2]. Now, the result follows from \cref{kan}. 
 \end{proof}

\medskip
The following theorem establishes a connection between closed and extension-closed subcalsses of $\nPSpec \cA$ and torsion classes of $\cA$.

\begin{Theorem}\label{tors}
There exists an injective map $\Phi\mapsto \nPSupp^{-1}\Phi$  between closed and extension-closed subclasses of $\nPSpec\cA$
and torsion classes of $\cA$. This map has a left inverse given by  $\cX\mapsto\nPSupp\cX$. 
\end{Theorem}
\begin{proof}
Let $\Phi$ be a closed and extension-closed subset of $\nPSpec\cA$ and let $\cX$ be a torsion class of $\cA$. We show that $\nSupp^{-1}\Phi$ is a torsion subclass of $\cA$ . Suppose that $\{M_i\}_{\Lambda}$ is  a family  of objects in $\nSupp^{-1}\Phi$. We show that 
$M=\coprod_{\Lambda} M_i\in\nSupp^{-1}\Phi$. Given $[H]\in\nPSupp M$, we have 
$[H]\in\overline{C(M)}$. Then $H$ has a finite filtration of subobjects  $0\subset H_1\subset\dots\subset H_{n-1}\subset H$ such that each $H_i/H_{i-1}$ is a quotient of $M$ for each $i$. Since $H$ is noetherian, there exists a finite subset $\Gamma$ of $\Lambda$ such that $H_i/H_{i-1}$ is a quotient of $\coprod_{\Gamma}M_j$ for each $i$. As $\Phi$ is extension-closed, $\coprod_{\Gamma}M_j\in\nPSupp^{-1}\Phi$. Thus $H_i/H_{i-1}\in\nPSupp^{-1}\Phi$ for each $i$ and consequently $H\in\nPSupp^{-1}\Phi$. Therefore $[H]\in\Phi$. On the other hand, it follows from [LS, Proposition 6.2] that $\nPSupp\cX=\nPSupp(\cX\cap\noeth\cA)$; and hence $\nPSupp\cX$ is closed and extension-closed subclass of $\nPSpec\cA$ by [LS, Theorem 6.8]. It is clear that $\nPSupp(\nPSupp^{-1}\Phi)=\Phi$. 
\end{proof}

\medskip

For a torsion class $\cT$ of $\cA$, a functor $t:\cA\To \cT$ is said to be {\it idempotent radical} if $t^2=t$ and  $t(X/t(X))=0$ for every object $X$ in $\cA$. According to [S, Chap VI, Proposition 2.3], there is a bijective correspondence between torsion subclasses of $\cA$ and idempotent radicals.  A torsion class $\cT$ of $\cA$ is said to be of {\it finite type} provided that the idempotent radical functor $t:\cA\To\cT$ commutes with direct limits. For every subcategory $\cS$ of $\cA$, we denote by $\overset{\To}\cS$, the full subcategory of $\cA$ consisting of $\underset{\rightarrow}{\rm lim}X_i$ with $X_i\in\cS$ for each $i$. A torsion class $\cT$ of $\cA$ is said to be of {\it weakly finite type} if $\cT=\overrightarrow{\cT\cap\noeth\cA}$. We now have the following lemma.

\medskip

\begin{Lemma}\label{wft}
Let $\cA$ be a locally noetherian Grothendieck category. Then every torsion class of finite type of $\cA$ is of weakly finite type.
\end{Lemma}
\begin{proof}
 Let $\cT$ be a torsion class of weakly finite type of $\cA$. Then the inclusion  $\overrightarrow{\cT\cap\noeth\cA}\subset \cT$ is clear. Now assume that $M\in\cT$ and so $M=\underset{\rightarrow}{\rm lim}M_i$ is the direct limit of its noetherian subobjects. Thus $M=t(M)\cong \underset{\rightarrow}{\rm lim}t(M_i)\in\overrightarrow{\cT\cap\noeth\cA} $ 
\end{proof}
\medskip
It is natural to ask whether over a locally noetherian Grothendieck category $\cA$, a torsion class of weakly finite type of $\cA$ is of finite type. However, this holds if $\cT$ is a localizing subcategory of $\cA$, see [Kr, Theorem 2.8].

\medskip 
\begin{Lemma}\label{wektor}
Let $\cA$ be a locally noetherian Grothendieck category and let $\cT$ be a torsion class of $\cA$. Then $\overrightarrow{\cT\cap\noeth A}$ is a torsion class of weakly finite type of $\cA$. 
\end{Lemma}
\begin{proof}
We first show that $\overrightarrow{\cT\cap\noeth A}$ is the torsion class of $\cA$ generated by $\cT\cap\noeth A$. By [S, Chap. VI, Proposition 2.5], it suffices to show that for every  object $M$ in $\overrightarrow{\cT\cap\noeth A}$, every nonzero quotient $M/N$ of $M$ has a nonzero subobject in $\cT\cap\noeth\cA$. Since $M\in\overrightarrow{\cT\cap\noeth A}$, we can write $M=\underset{\rightarrow}{\rm lim}M_i$ where $M_i\in\cT\cap\noeth\cA$ for each $i$. Suppose that $\{\alpha_i|\hspace{0.1cm} \alpha_i:M_i\To M\}$  is the canonical compatible family of morphisms and $\pi:M\To M/N$ is the canonical epimorphism. If for each $i$, we have $\pi\alpha_i=0$, then $\Im\alpha_i\subset N$. Thus  $M=\Sigma_i\Im\alpha_i\subset N$ which is a contradiction. Therefore there exists some $i$ such that $\pi\alpha_i\neq 0$; and hence $M/N$ contains $\Im \pi\alpha_i\in\cT\cap\noeth\cA$. To prove that $\overrightarrow{\cT\cap\noeth A}$ is a weakly finite type of $\cA$, it suffices to show that $\overrightarrow{\cT\cap\noeth A}\cap\noeth\cA=\overrightarrow{\cT\cap\noeth A}$. If $M\in\overrightarrow{\cT\cap\noeth A}\cap \noeth \cA$, it follows from [B, Theorem 4.1] that $M$ is a direct summand of some object in $\cT\cap\noeth\cA$ so that $M\in\cT\cap\noeth\cA$.
\end{proof}

 The following theorem classifies torsion classes of weakly finite type of a locally noetherian Grothendieck category. 
 
\medskip

\begin{Theorem}\label{classloc}
 Let $\cA$ be a locally noetherian Grothendieck category. Then the map $\cX\mapsto\nPSupp\cX$ establishes a one-to-one correspondence between torsion classes of weakly finite type of $\cA$ and closed and extension-closed subsets of $\nPSpec\cA$. The inverse map is given by $\Phi\mapsto \overrightarrow{\nPSupp^{-1}\Phi\cap\noeth A}$.
\end{Theorem}
\begin{proof}
Suppose that $\cX$ is a torsion class of weakly finite type of $\cA$ and  suppose that $\Phi$ is a closed and extension-closed subset of $\nPSpec \cA$. According to \cref{tors}, $\nPSupp\cX$ is a closed and extension-closed subset of $\nPSpec\cA$ and $\nPSupp^{-1}\Phi$ is a torsion class of $\cA$. Then $ \overrightarrow{\nPSupp^{-1}\Phi\cap\noeth A}$ is a torsion class of weakly finite type of $\cA$ by \cref{wektor}. We now show that $\nSupp(\overrightarrow{\nPSupp^{-1}\Phi\cap\noeth \cA})=\Phi$. The inclusion $\Phi\subset\nSupp(\overrightarrow{\nPSupp^{-1}\Phi\cap\noeth A})$ is clear. For other side, if $[H]\in\nSupp(\overrightarrow{\nPSupp^{-1}\Phi\cap\noeth A})$, then there exists an object $M\in\overrightarrow{\nPSupp^{-1}\Phi\cap\noeth\cA}$ such that $[H]\in\nPSupp M$. Therefore $M\in\nPSupp^{-1}\Phi$ as $\nPSupp^{-1}\Phi$ is a torsion class of $\cA$. Consequently $[H]\in\Phi$. We also show that $
\overrightarrow{\nPSupp^{-1}(\nPSupp\cX)\cap\noeth\cA}=\cX$. If $M\in
\overrightarrow{\nPSupp^{-1}(\nPSupp\cX)\cap\noeth\cA}$, then $M=\underset{\rightarrow}{\rm lim}M_i$, where $M_i\in \nPSupp^{-1}(\nPSupp\cX)\cap\noeth\cA$ for each $i$. Thus $\nPSupp M_i\subset \nPSupp \cX$ for each $i$. It follows from [LS, Proposition 6.2] that $M_i\in\cX$ for each $i$ and so $M\in\cX$ as $\cX$ is a torsion class of $\cA$. For other side, given $M\in\cX$, since $\cX$ is weakly finite type of $\cA$, we have $M=\underset{\rightarrow}{\rm lim} M_i$, where $M_i\in\cX\cap\noeth\cA$ for each $i$.  Thus  $M_i\in\nPSupp^{-1}(\nPSupp\cX)\cap\noeth \cA$ and so $M\in \overrightarrow{\nPSupp^{-1}(\nPSupp\cX)\cap\noeth\cA}$.
\end{proof}

\medskip

In the rest of this section, we assume that the collection of all closed and extension-closed subclasses of $\nPSpec\cA$ forms a topology on $\nPSpec\cA$. We have the following lemma.

\begin{Lemma}\label{extlem}
 Let $0\To X\To Y\To Z\To 0$ be an exact sequence of noetherian objects in $\cA$. Then $\nSupp Y\subset \nPSupp X\cup\nPSupp Z$. In particular, every closed subclass of $\nPSpec \cA$ is extension-closed. 
   \end{Lemma}
   \begin{proof}
 By \cref{exttt} and the assumption $\nPSupp X\cup\nPSupp Z$ is an extension-closed subclass of $\nPSpec\cA$ and so $\nSupp Y\subset \nPSupp X\cup\nPSupp Z$. The second claim is clear using the first part.
\end{proof}

\medskip

 A closed subclass $\Phi$ of $\PSpec\cA$ is {\it compact} if every closed cover of $\Phi$ admits a finite subcover. Under a certain condition, the closed compact subclasses of $\PSpec\cA$ can be characterized as follows.
 
\begin{Proposition}\label{compp}
Let the collection of all closed and extension-closed subclasses of $\nPSpec\cA$ form a topology on $\nPSpec\cA$. Then a closed subclass $\Phi$ of $\nPSpec\cA$ is compact if and only if there exists a noetherian object $M$ in $\cA$ such that $\Phi=\nPSupp M$.
\end{Proposition}
\begin{proof}
Assume that $M$ is a noetherian object in $\cA$ and assume that $\{\Phi_{\lambda}\}_{\lambda\in\Lambda}$ is a cover of $\PSupp M$. According to \cref{filt}, there exists a finite filtration of subobjects $$0=M_0\subset M_1\subset \dots \subset M_{n-1}\subset M_n=M$$ of $M$ and noetherian premonoform objects $H_1,\dots H_n$ in $\cA$ such that $[H_1],\dots[H_n]\in\nPSupp M$ and  $M_i/M_{i-1}$ is a quotient of $H_i$ for $i=1,\dots,n$. Since $\PSupp M=\bigcup_{\lambda\in\Lambda}\Phi_{\lambda}$, there exists ${\lambda_1},\dots, \lambda_n\in\Lambda$ such that $\PSupp H_i\subset \Phi_{\lambda_i}$ for $i=1,\dots,n$. By virtue of \cref{exttt}, each $\nPSupp H_i$ is closed and extension-closed and so the assumption implies that $\bigcup_{i=1}^n\nPSupp H_i$ is an closed and extension-closed subset of $\nPSpec\cA$.  Hence $\nPSupp M\subset\bigcup_{i=1}^n\nPSupp H_i\subset \bigcup_{i=1}^n\Phi_{\lambda_i}.$ Conversely assume that $\Phi\subset\nPSpec\cA$ is compact. Clearly $\Phi=\bigcup_{[G]\in\Phi}\nPSupp G$ and so there exists $[G_1],\dots [G_t]\in\Phi$ such that $\Phi= \bigcup_{i=1}^t\PSupp G_i$. We observe that $G_i$ is noetherian for $i=1,\dots,n$. Using the assumption, $\Phi$ is an extension-closed subset of $\nPSpec\cA$ so that   $\Phi=\nPSupp\bigoplus_{i=1}^nG_i$.
\end{proof}

\section{Premonoform modules over a noncommutative ring}

In this section, we assume that $A$ is a ring with identity.
A two-sided ideal $\frak p$ of $A$ is {\it extremely  prime} if for every $a,b$ in $A$, the condition $ab\in\frak p$ implies that $a\in\frak p$ or $b\in\frak p$. (Such ideals have been studied, for instance in [A]). For an $A$-bimodule $M$ and $x\in M$,  the {\it left annihilator} of $x$ is $\lAnn_A x=\{r\in A|\hspace{0.1cm} rx=0\}$ and the left annihilator of $M$ is $\lAnn_AM=\{r\in A|\hspace{0.1cm} rM=0\}$. It is clear that $\lAnn_AM$ is a two-sided ideal. The {\it right annihilator} of $x$ ($M$ resp.), denoted by ${\rAnn_Ax}$ ($\rAnn_AM$ resp.) is defined similarly.

\medskip

\begin{Proposition}\label{prmonex}
Let $\frak p$ be a two-sided ideal of $A$. Then the following statements are equivalent.

${\rm (i)}$ $\frak p$ is extremely prime.

${\rm (ii)}$ $A/\frak p$ is a monoform right $A$-module. 
 
${\rm (iii)}$ $A/\frak p$ is a premonoform right $A$-module.
 
 ${\rm (iv)}$ $A/\frak p$ is a monoform left $A$-module.
 
 ${\rm (v)}$ $A/\frak p$ is a premonoform left $A$-module.
\end{Proposition}
\begin{proof}
We prove (i)$\Rightarrow$ (ii)$\Rightarrow$ (iii)$\Rightarrow$ (i) and  (i)$\Rightarrow$ (iv)$\Rightarrow$ (v)$\Rightarrow$ (i) is similar. To prove (i)$\Rightarrow$ (ii) if $A/\frak p$ is not monoform, then there exists  a right ideal $J$ of $A$ containing $\frak p\subsetneq J$ such that $A/\frak p$ has a nonzero submdule $N$ isomorphic to a submodule $K$ of  $A/J$. Then for every nonzero element $a+J$ of $K$, we have $J\subset \lAnn_A(a+J)=\frak p$  which is a contradiction. (ii)$\Rightarrow$ (iii) is clear.  (iii)$\Rightarrow$ (i) Assume that  $a,b\in A$ and $a,b\notin\frak p$. We assert that $\frak p=\rAnn_R(a+\frak p)$. The inclusion $\frak p\subset \rAnn_A(a+\frak p)$ is clear. Assume that $x\in\rAnn_A(a+\frak p)$. Since $a\notin\frak p$, the homomorphism $A/\frak p\stackrel{a.}\to A/\frak p$ of right $A$-modules is injective because $A/\frak p$ is premonoform. This implies that $x\in\frak p$. Thus $\frak p=\rAnn_A(a+\frak p)\subsetneq \rAnn_A(1+bA+\frak p)$ as $A/bA+\frak p$ is not isomorphic to a submodule of $A/\frak p$. Then there exists $y\in \rAnn_A(1+bA+\frak p)$ such that $ay\notin\frak p$. Moreover, there exist  $r\in A$ and $p\in\frak p$ such that $y=br+p$. Hence the equality $ay=abr+ap$ implies that $ab\notin\frak p$.
\end{proof}

\medskip
\begin{Proposition}\label{pp}
Suppose that $M$ is an $A$-bimodule that is a premonoform right $A$-module. Then

${\rm (i)}$ $\lAnn_A x=\lAnn_AM$ for every nonzero element $x\in A$. 

 ${\rm (ii)}$ $\lAnn_AM$ is an extremely prime ideal of $A$.
 
 ${\rm (iii)}$ If $A$ is left noetherain and $M$ is a uniform left $A$-module, then $M$ is a monoform left $A$-module.
 \end{Proposition} 
\begin{proof}
(i) Assume that $x\in M$ is a nonzero element and  $r\in\lAnn_Ax$ such that $r\notin\lAnn_AM$.  Since $M$ is a premonoform right $A$-module, the endomorphism $M\stackrel{r.}\To M$ of right $A$-module $M$ is injective. Thus $x=0$ which is a contradiction. (ii) Assume that $rs\in\lAnn_RM$ and $r\notin\lAnn_AM$. Then the endomorphism $M\stackrel{r.}\To M$ of right $A$-module $M$ is injective and so $s\in\lAnn_AM$. (iii) We first show that every noetherian submodule $N$ of $M$ is monoform.  Using an easy induction, we may assume that $N=Ax+Ay$. We notice that $Ax\cong A/\lAnn_Ax$ and $Ay\cong A/\lAnn_Ay$. It follows from (i), (ii) and \cref{prmonex} that $Ax$ and $Ay$ are monoform left $A$-modules. Thus $N$ is a monoform left $A$-module by [K, Proposition 2.11]. In order to prove that $M$ is a monoform left $A$-module, assume that $M$ has a nonzero submodule $D$ such that  $M$ has a nonzero submodule $Ax$ and $M/D$ has a nonzero submodule $K/D$ with $K$ noetherian and $Ax\cong K/D$. Since $M$ is uniform, $K\cap Ax$ is a common  nonzero subject of $K$ and $K/D$. But this implies that $K$ is not monoform which is a contradiction. 
\end{proof}

\begin{Proposition}\label{pp1}
Suppose that $M$ is an $A$-bimodule that is a premonoform left $A$-module. Then

${\rm (i)}$ $\rAnn_A x=\rAnn_AM$ for every nonzero element $x\in A$. 

 ${\rm (ii)}$ $\rAnn_AM$ is an extremely prime ideal of $A$.
 
 ${\rm (iii)}$ If $A$ is a right noetherian and $M$ is a uniform right $A$-module, then $M$ is a monoform right $A$-module.
 \end{Proposition} 
\begin{proof}
The proof is similar to that of \cref{pp}. 
\end{proof}

\medskip
\begin{Corollary}\label{bipr}
Suppose that $A$ is right and left noetherian ring and $M$ is an $A$-bimodule. If $M$ is a uniform and  premonoform left and right $A$-module, then $M$ is a monoform left and right $A$-module.
\end{Corollary}
\begin{proof}
Straightforward by \cref{pp} and \cref{pp1}.
\end{proof}

\begin{Example}
Let $k$ be a field a and $V$ be a vector space over $k$ with the basis $\{a,b\}$. Consider the tensor algebra $T(V)=\bigoplus_{i=0}^{\infty}V^{\otimes i}$ of  $V$. We observe that $T(V)$ is a $T(V)$-bimodule. According to [LS, Proposition 4.3] and a dual result, $T(V)$ is a premonoform left $T(V)$-module and a premonoform right $T(V)$-module but it is not a monoform left $T(V)$-module nor a monoform right $T(V)$-module. Furthermore, \cref{bipr} implies that $T(V)$ is neither a uniform left $T(V)$-module nor a uniform $T(V)$ right $T(V)$-module.
\end{Example}

Liu and Stanely [LS, Proposition 4.7] proved that over a commutative noetherian ring, every noetherian premonoform module is monoform. A parall result can be given as follows.

\medskip
\begin{Corollary}
Suppose that $A$ is a commutative noetherian ring and $M$ is a uniform and premonoform  $A$-module, then $M$ is  monoform.
\end{Corollary}

We show that completely prime right ideals in [R] are in correspondence to cyclic premonoform right $A$-modules. We recall that a right ideal $\frak a$ of $A$ is a {\it completely prime right ideal} of $A$ if for any $x,y\in A$ the conditions $x\frak a\subset\frak a$  and $xy\in\frak a$ imply $x\in\frak a$ or $y\in\frak a$.

\begin{Proposition}\label{compr}
Let $\frak a$ be a right ideal of $A$. Then $A/\frak a$ is a premonoform  right $A$-module if and only if $\frak a$ is  completely prime.
\end{Proposition}
\begin{proof}
Assume that $x,y\in A$ such that $x\frak a\subset\frak a$ and $xy\in\frak a$. Then $A/\frak a\stackrel{x.}\To A/\frak a$ is a homomorphism of right $A$-modules. Since $A/\frak a$ is premonoform, either $x.=0$ or $\ker x.=0$ so that $x\in\frak a$ or $y\in\frak a$. Conversely, for every nonzero homomorphism $f\in\End_A(A/\frak a)$, we have $f=x.$ where $f(1+\frak a)=x+\frak a$ and $x\notin\frak a$. Clearly $x\frak a\subset \frak a$ and since $\frak a$ is completely prime, $f$ is injective. 
\end{proof}

 \section{A  new topological space for an abelian category}
 Let $\cA$ be an abelian category. In this section we define a new topology on $\PSpec\cA$ in which the closed subsets is defined using Serre subcategories.
 \begin{Definition}
For every object $M$ of $\cA$, we define the {\it Serre psupport} of $M$ as 
 $$\sPSupp M=\{[H]\in\PSpec\cA \hspace{0.05cm}|\hspace{0.15cm} H\in\sqrt{M}\}.$$ It is clear by the definition that $\PSupp M\subset\sPSupp M$. We also define {\it noetherian Serre support} of $M$ as  $\nsPSupp M=\sPSupp M\cap \nPSpec\cA$. If $M$ is a noetherian object in $\cA$, then $\nsPSupp M=\sPSupp M$. A subclass $\Phi$ of $\PSpec\cA$ is said to be {\it Serre closed} if for every $[H]\in\Phi$, we have $\sPsupp H\subset\Phi$. The subclass $\Phi$ is {\it Serre extension-closed} if for every exact sequence $0\To M\To X\To N\To 0$ in $\cA$, the condition 
 $\sPSupp M,\sPSupp N\subset \Phi$ implies that $\sPSupp X\subset \Phi$.  
 \end{Definition}
 \medskip
 
 \begin{Lemma}\label{opp}
 Let $M$ be an object in $\cA$. Then

 ${\rm (i)}$ $\sPSupp M$ is a Serre closed subclass of $\cA$. 
 
  ${\rm (ii)}$ If $N$ is a subobject of $M$, then $\sPSupp N,\sPSupp M/N\subset \sPSupp M$.
  \end{Lemma}
  \begin{proof}
  (i) If $[H]\in\sPSupp M$, then $H\in\sqrt{M}$. Hence $\sqrt{H}\subset \sqrt{M}$. It is straightforward by definition that $\sPSupp H\subset\sPSupp M$. (ii) If $[H]\in\sPSupp N$, then $H\in\sqrt{N}\subset \sqrt{M}$. Hence $[H]\in\sPSupp M$. The proof for $\sPSupp M/N$ is similar.
  \end{proof}
 \medskip
 For any subcategory $\cX$, the Serre psupport of $\cX$ is defined as $\sPSupp\cX=\bigcup_{M\in\cX}\sPSupp M$. Moreover $\nsPSupp\cX=\sPSupp\cX\cap\nPSpec\cA$. We observe that if $\cX\subset\noeth \cA$, then $\nsPSupp\cX=\sPSupp\cX$.

 \medskip
 \begin{Example}
 Let $p$ be a prime number. For every $[H]\in\nsPSupp\bbZ_{p^{\infty}}$, the abelian group $H$ is  monoform and it contains a subobject $H_1$ which is a subqoutient of $\bbZ_{p^{\infty}}$. Since every qoutient of $\bbZ_{p^{\infty}}$ is injective, we deduce that $\overline{H}=\overline {\bbZ/p\bbZ}$. Thus \cref{isosub} implies that 
 $\nsPSupp\bbZ_{p^{\infty}}=\{[\bbZ/p\bbZ]\}$.  
 \end{Example}

 \medskip
 
 \begin{Lemma}\label{ssur}
 Let $\cX$ be a subcategory of $\cA$ closed under finite direct sums. Then $$\sPSupp \cX=\{[H]\in\PSpec \cA|\hspace{0.1cm} H\in\sqrt{\cX}\}.$$
 \end{Lemma}
 \begin{proof}
 The inclusion $\sPSupp \cX\subset\{[H]\in\PSpec \cA|\hspace{0.1cm} H\in\sqrt{\cX}\}$ is clear. Assume that $[H]\in\PSpec\cA$ such that $H\in\sqrt{\cX}$. Then $H$ has a filtration $0\subset H_1\subset\dots \subset H_{n-1}\subset H_n=H$ and $X_i\in\cX$ such that $H_i/H_{i-1}$ is a subquotient of $X_i$ for $i=1,\dots,n$. By the assumption, we have $X=\bigoplus_{i=1}^nX_i\in\cX$ and $H\in\sqrt{X}$ so that $[H]\in\sPSupp X$. 
 \end{proof}
 \medskip

 \begin{Corollary}\label{corcom}
 Let $A$ be a commutative noetherian ring and $\cX$ be a subcategory of noetherian $A$-modules closed under finite direct sums. Then $\nsPSupp\cX=\nPSupp\cX$. In particular, $\nsPSupp M=\nPSupp M$ for every noetherian $A$-module $M$.
 \end{Corollary}
 \begin{proof}
 It follows from \cref{sertor} that $\overline{C(\cX)}=\sqrt{\cX}$. Now the result is obtained using \cref{sursum} and \cref{ssur}.  To prove the second assertion, using the first part, we have $\nPSupp M=\nPSupp \overline{C(M)}={\sqrt{M}}=\nsPSupp M$.
 \end{proof}

 \medskip
 \begin{Proposition}\label{topo}
 The collection of the Serre closed subclasses of $\PSpec \cA$ forms an Alexandroff topology.
 \end{Proposition}
 \begin{proof}
 If $\{\Phi_i\}_I$ is a family of Serre closed subclasses of $\PSpec\cA$, it is straightforward by definition that $\bigcup_I\Phi_i$ and $\bigcap_I\Phi_i$ are Serre closed subclasses of $\PSpec\cA$.
 \end{proof}

 \medskip
 
 \begin{Corollary}\label{comclex}
 Let $A$ be a commutative noetherian ring. Then closed subsets and Serre closed subsets of $\nPSpec A$ are the same. Furthermore, extension-closed subsets and Serre extension subsets of $\nPSpec A$ are the same.
 \end{Corollary}
 \begin{proof}
 Let $\Phi$ be a closed subset of of $\nPSpec A$. Then $\Phi=\bigcup_{[H]\in\Phi}\nPSupp H=\bigcup_{[H]\in\Phi}\nsPSupp H$ by \cref{corcom}. Therefore $\Phi$ is a Serre closed subset of $\nPSpec A$ by  \cref{topo}. A similar argument deduces that a Serre closed subset of $\nPSpec A$ is a closed subset of $\nPSpec A$. The proof of the second claim is straightforward using the first part.  
 \end{proof}
 
 \medskip

  \begin{Lemma}\label{Ser}
 Let $\cA$ be an abelian category, let $\cS$ be a Serre subcategory of $\cA$ and let $M$ be a noetherian object in $\cA$. If $\sPSupp M\subset \sPSupp\cS$, then $M\in\cS$. 
 \end{Lemma}
 \begin{proof}
 Since $M$ is noetherian, by virtue of \cref{filt}, it has a finite filtration of submodules
$$0=M_0\subset M_1\subset \dots \subset M_{n-1}\subset M_n=M$$ such that $M_i/M_{i-1}$ is a quotient of some noetherian premonoform object $H_i\in\cA$ with $[H_i]\in\PSupp M$ for each $i$. Then  $[H_i]\in\sPSupp M$ and the assumption implies that $H_i\in\sqrt{X_i}$ for some $X_i\in\cS$ for each $i$. Therefore $H_i\in\cS$ for each $i$; and consequently $M\in\cS$.
 \end{proof}
 \medskip
 
 \begin{Lemma}\label{Serrr}
 Let $\cS$ be a Serre subcategory of $\noeth\cA$. Then $\nsPSupp\cS$ is a Serre extension-closed subclass of $\nPSpec\cA$. In particular, $\nsPSupp M$ is a Serre extension-closed subclass of $\nPSpec \cA$ for every noetherian object $M$ in $\cA$. 
 \end{Lemma}
 \begin{proof}
 Assume that $0\To M\To X\To N \To 0$ is an exact sequence of objects in $\noeth\cA$ such that $\sPSupp M,\sPSupp N\subset\sPSupp \cS$. It follows from \cref{Ser} that $M,N\in\cS$ so that $X\in\cS$. Hence $\sPSupp X\subset\sPSupp\cS$. The second claim follows from the first part using the fact that $\nsPSupp M=\nsPSupp\sqrt{M}$. 
 \end{proof}
 
 \medskip
  \begin{Theorem}
The map $\cS\mapsto \nsPSupp \cS$ establishes a one-to-one correspondence between Serre subcategories of $\noeth\cA$ and Serre closed and Serre extension-closed subclasses of $\nPSpec\cA$. The inverse map is given by $\Phi\mapsto \nsPSupp^{-1}\Phi$. 
  \end{Theorem}
  \begin{proof}
  It is straightforward using \cref{Ser} and \cref{Serrr}. 
  \end{proof}
 
 \medskip
 In the rest of this section, we assume that $\cA$ is a locally noetherian Grothendieck category. The following theorem classifies the localizing subcategoris of $\cA$ via the Serre closed and Serre extension closed subsets of $\nPSpec\cA$.

 \medskip
 
 \begin{Theorem}
 The map $\cX\mapsto\nsPSupp\cX$ establishes a one-to-one correspondence between localizing subcategories of $\cA$ and Serre closed and Serre extension-closed subsets of $\nPSpec\cA$.  The inverse map is given by $\Phi\mapsto\nsPSupp^{-1}\Phi$.
 \end{Theorem}
 \begin{proof}
 Assume that $\cX$ is a localizing subcategory of $\cA$. It is clear by \cref{Ser} that $\nsPSupp\cX=\nsPSupp(\cX\cap\noeth\cA)$. Hence $\nsPSupp\cX$ is a is Serre extension-closed subset of $\nPSpec\cA$ by \cref{Serrr}. now, assume that $\Phi $ is a Serre closed and Serre extension-closed subset of $\nPSpec\cA$. Clearly $\nsSupp^{-1}\Phi$ is a Serre subcategory of $\cA$. If $\{M_i\}_I$ is an arbitrary family of objects in $\nsSupp^{-1}\Phi$ and $[H]\in\nsSupp\coprod_I M_i$, then by \cref{kan}, there is a finite filtration $0=H_0\subset H_1\subset H_1\subset\dots\subset H_n=H$ of subobjects in $H$ such that each $H_j/H_{j-1}$ is a subquotient of $\coprod_I M_i$. Since each $H_j/H_{j-1}$ is noetherian, there exists a finite subset $J$ such that $H_j/H_{j-1}$ is a subquotient of $M_J=\coprod_J M_i$  for each $j$. As $\Phi$ is Serre extension-closed, $M_J\in\nsPSupp^{-1}\Phi$. Thus $H_j/H_{j-1}\in\nsSupp^{-1}\Phi$ for each $j$ which forces that $H\in\nsSupp^{-1}\Phi$. Then $[H]\in \Phi$; and hence $\coprod M_i\in\nsSupp^{-1}\Phi$. We now show that $\nsSupp^{-1}(\nsSupp\cX)=\cX$. Assume that $M\in\nsSupp^{-1}(\nsSupp\cX)$. Since $\cA$ is locally noetherain, $M=\bigcup_i M_i$ is the direct union  of its noetherian subobjects. Then $M_i\in\nsSupp^{-1}(\nsSupp\cX)$ for each $i$ so that $M_i\in\cX$ using \cref{Ser}. Thus $M\in\cX$ and  so $\nsSupp^{-1}(\nsSupp\cX)\subset\cX$. The other side is clear. The equality $\nsPSupp(\nsPSupp^{-1}\Phi)=\Phi$ is straightforward.
 \end{proof}
 
 \medskip
 
 \begin{Proposition}\label{exten}
 Let $\cX$ be a localizing subcategory of $\cA$. Then $\sPSupp\cX$ is a Serre extension-closed subset of $\PSpec\cA$.
 \end{Proposition}
 \begin{proof}
 Let $0\To N\To M\To K\To 0$ be an exact sequence of objects of $\cA$ such that $\sPSupp N,\sPSupp K\subset\sPSupp\cX$. Since $\cA$ is locally  noetherian, $M=\bigcup _iM_i$ is the direct union of its noetherian subobjects. For each $i$, there is an exact sequence $0\To N_i\To M_i\To K_i$ such that $N_i$ is a subobjects of $N$ and $K_i$ is a subobject of $K$. Also $N=\bigcup_iN_i$ and $K=\bigcup_iK_i$.
 The assumption implies that $\sPSupp N_i,\sPSupp K_i\subset\sPSupp\cX$ and so  \cref{Ser} implies that $N_i,K_i\in\cX$ for each $i$. Thus $N,K\in\cX$ so that $M\in\cX$. Consequently, $\sPSupp M\subset\sPSupp\cX$.
   \end{proof}
 
  In the rest of this section, we assume that  the collection of Serre closed and Serre  extension-closed subsets as the closed subsets of $\nPSpec\cA$ forms a topology on $\nPSpec\cA$. For example, if $A$ is a commutative noetherian ring, then $\nPSpec A$ satisfies this condition by \cref{extclos} and \cref{comclex}.  
  
 \medskip
 
 \begin{Proposition}\label{ext}
 Let $0\To X\To Y\To Z\To 0$ be an exact sequence of objects in $\cA$. Then $\nsSupp Y=\nsPSupp X\cup\nsPSupp Z$. In particular,  every Serre closed subset of $\nPSpec \cA$ is Serre extension-closed. 
   \end{Proposition}
   \begin{proof}
 We first assume that $Y$ is noetherian. Then by \cref{Serrr} and the assumption in the beginning of the proposition, $\nPSupp X\cup\nsPSupp Z$ is a Serre extension-closed subset of  $\nPSpec\cA$ so that $\nsPSupp Y\subset \nsPSupp X\bigcup\nsPSupp Z$. It now follows from \cref{opp} that $\nsPSupp Y=\nsPSupp X\bigcup\nsPSupp Z$. If $Y$ is an arbitrary object, then $Y=\bigcup_IY_i$ is the direct union of its noetherian subobjects $Y_i$. For each $i$, there is an exact sequence $0\To X_i\To Y_i\To Z_i\To 0$ such that $X_i$ is a subobject of $X$ and $Z_i$ is a subobject of $Z$ and $X=\bigcup_IX_i$ and $Z=\bigcup_IZ_i$. According to the first case, we have $\nsPSupp Y_i=\nsPSupp X_i\cup \nsPSupp Z_i$ for each $i$. Given $[H]\in\nsPSupp Y$, the noetherian object $H$ has a finite filtration of its subobjects $H$$$0=H_0\subset H_1\subset\dots \subset H_n=H$$  such that each $H_j/H_{j-1}$ is a subquotient of $Y$. Since each $H_j/H_{j-1}$ is noetherian, there exists $\alpha\in I$ such that $H_j/H_{j-1}$ is a subquotient of $Y_{\alpha}$. Finally by the assumption and the firs case, we have $\nsPSupp H\subset\nsPSupp Y_{\alpha}=\nsPSupp X_{\alpha}\cup\nsPSupp Z_{\alpha}\subset \nsPSupp X\cup\nsPSupp Z$ so that $[H]\in\nsPSupp X\cup\nsPSupp Z$. Thus $\nsPSupp Y\subset\nsPSupp X\cup\nsPSupp Z$ and so the equality holds by \cref{opp}. The second assertion is clear by the first part.     
\end{proof}

\medskip
\begin{Proposition}
Let $\Phi$ be a Serre closed subset of $\nPSpec\cA$. Then $\Phi$ is Serre extension-closed and compact if and only if there exists a noetherian object $M$ in $\cA$ such that $\Phi=\nsPSupp M$. 
\end{Proposition}
\begin{proof}
The proof is similar to the proof of \cref{compp}.
\end{proof}
 
 \section{Premonoform modules over a commutative ring} 

In this section, we assume that $A$ is a commutative ring with identity, $\Mod A$  denotes the category of $A$-modules and mod-$A$ denotes the subcategory of finitely generated $A$-modules. Furthermore, we denote $\PSpec \Mod A$ (resp. $\nPSpec\Mod A$) by $\PSpec A$ (resp. $\nPSpec A$). For every $A$-module $M$, we set $\frak p_M=\Ann_AM$. According to \cref{pp,pp1}, $\frak p_H$ is a prime ideal of $A$ for every premonoform $A$-module $H$.
\medskip
\begin{Proposition}\label{f}
 The following conditions hold.

${\rm(i)}$ For any prime ideals $\frak p$ and $\frak q$ of $A$, we have $[A/\frak p]=[A/\frak q]$ if and only if $\frak p=\frak q$. 

${\rm(ii)}$ For every ideal $I$ of $A$,  $\PSupp  A/I\subset\{[H]\in\PSpec A|\hspace{0.1cm} \frak p_H\in V(I)\}$. If $A$ is noetherian, then $\PSupp  A/I=\{[H]\in\nPSpec A|\hspace{0.1cm} \frak p_H\in V(I)\}$.

${\rm(iii)}$ For every finitely generated $A$-module $M$, $\PSupp M\subset\{[H]\in\PSpec A|\hspace{0.1cm} \frak p _{H}\in\Supp M\}$. If $A$ is noetherian, then $\PSupp M=\{[H]\in\nPSpec A|\hspace{0.1cm} \frak p _{H}\in\Supp M\}$. 
\end{Proposition}
\begin{proof}
(i) If $[A/\frak p]=[A/\frak q]$, then $\overline{C(A/\frak p)}=\overline{C(A/\frak q)}$. Since $A/\frak p\in\overline{C(A/\frak q)}$, it has a finite filtration of submodules $$0=X_0\subset X_1\subset\dots\subset X_n=A/\frak p$$ such that $X_i/X_{i-1}$ is a quotient of some module isomorphic to $A/\frak q$ for $i=1,\dots,n$. In particular $X_1\cong A/J$ such that $\frak q\subset J$. Since $X_1$ is monoform, $J=\Ann X_1=\frak p$ so that $\frak q\subset \frak p$. A similar argument deduces that $\frak q\subset \frak p$. (ii)  Assume that $[H]\in\PSupp A/I$. Similar to the proof of (i), the module $H $ has a submodule $X$ which is a quotient of some module isomorphic to $A/I$. Thus $I\subset\frak p_H=\Ann_AX$. If $A$ is noetherian, then $\PSupp A/I \subset \nPSpec\cA$ and so $\PSupp  A/I\subset\{[H]\in\nPSpec A|\hspace{0.1cm} \frak p_H\in V(I)\}$. To prove the other side, if $[H]\in\nPSpec A$ such that $I\subset \frak p_H$, then $H$ is a noetherian $A/I$-module.  Thus $H\in\overline{C(A/I)}$ so that $[H]\in\PSupp A/I$.  
(iii) is proved similar to (ii). 
\end{proof}

\begin{Lemma}\label{ann}
If $M,N$ are premonoform $A$-modules such that $[M]=[N]$, then $\frak p_M=\frak p_N$. 
\end{Lemma}
\begin{proof}
 Since $M\in\overline{C(N)}$, it admits a finite filtration of right submodules $$0=M_0\subset M_1\subset \dots\subset M_{n-1}\subset M_n=M$$
 such that $M_i/M_{i-1}$ is a quotient of some premonoform module $N_i$ isomorphic to $N$ for  $i=1,\dots,n$. Then for every $r\in\frak p_N$, it is easy to show that $r^nM=0$ and so $r\in\frak p_ M$. Consequently $\frak p_N\subset\frak p_M$. The other side is obtained similarly.  
\end{proof}

\medskip

\begin{Proposition}\label{com}
Let $M$ and $N$ be monoform $A$-modules. If $[M]=[N]$, then $\overline{M}=\overline{N}$. Furthermore, the converse  holds if $M$ and $N$ are noetherian,   
\end{Proposition}
\begin{proof}
As $M,N$ are monoform, there exist prime ideals $\frak p,\frak q$ of $A$ such that $\Ass_A M=\{\frak p\}$ and $\Ass_R N=\{\frak q\}$. It suffices to show that $\frak p=\frak q$. It follows from \cref{pp} that $\frak p=\frak p_M$  and $\frak q=\frak q_N$. Since $[M]=[N]$, by \cref{ann}, we have $\frak p_N=\frak p_M$. To prove the converse, since $\overline{M}=\overline{N}$, there exists a prime ideal $\frak p$ of $A$ such that $\Ass_AM=\Ass_AN=\{\frak p\}$. It suffices to prove that $[M]=[A/\frak p]$ and so a similar argument holds for $N$. We notice that $M$ is a faithful $A/\frak p$-module. Thus by  Gruson theorem [V, Theorem 4.1, p.82], the $A/\frak p$-module $A/\frak p$ admits  a finite filtration of submodules $$0\subset L_1\subset \dots\subset L_{n-1}\subset L_n=A/\frak p$$ such that each $L_i/L_{i-1}$ is a quotient of finite direct sum of $M$. Therefore $A/\frak p \in\overline{C(M)}$ and hence $\overline{C(A/\frak p)}\subset\overline{C(M)}$. Conversely, since $M$ is a noetherian $A/\frak p$-module, $M\in\overline{C(A/\frak p)}$ and so $\overline{C(M)}\subset \overline{C(A/\frak p)}.$  
\end{proof}
We recall from [Sm] that a nonzero $A$-module $M$ of $\cA$ is called {\it compressible} if each nonzero submodule $L$ of $M$ has some submodule isomorphic to $M$. 
\medskip
\begin{Corollary}
Let $M$ and $N$ be noetherian monoform $A$-modules such that $\overline{M}=\overline{N}$. Then $M$ is isomorphic to a submodule of $N$.  In particular, every noetherian monoform module is compressible. 
\end{Corollary}
\begin{proof}
The result is obtained using \cref{com} and \cref{isosub}. To prove the second claim, if $K$ is a submodule of $M$, then $\overline{K}=\overline{M}$. Thus using the first part, $K$ has a submodule isomorphic to $M$.   
\end{proof}

A subset $V$ of $\Spec A$ is said to be a {\it specialization-closed subsets} of $\Spec A$ if for prime ideal $\frak p,\frak a$ in $\Spec A$ with $\frak p\subset \frak q$, the condition $\frak p\in V$ yields that $\frak q\in V$. The following result establish a relationship between the specialization-closed subsets of $\Spec\cA$ and the closed and extension-closed subset of $\nPSpec \cA$.
\medskip
\begin{Proposition}\label{spec}
 Assume that $A$ is noetherian, $\Phi\subset \nPSpec A$ and
$\Psi=\{\frak p_H\in\nPSpec A|\hspace{0.1cm}[H]\in\Phi\}$. Then $\Psi$ is an specialization-closed subset of $\Spec A$ if and only if $\Phi$ is a closed and extension-closed subset of $\nPSpec A$.
\end{Proposition}
\begin{proof}
 Assume that $[H]\in\Phi$. We show that $\PSupp H\subset\Phi$. Given $[G]\in\PSupp H$, we have $G\in\overline{C(H)}$. Then $G$ has a submodule $G_1$ which is a quotient of $H$. Thus $\frak p_H\subset \frak p_{G_1}=\frak p_G$. Since $\frak p_H\in\Psi$, we deduce that $\frak p_{G}\in\Psi$. Hence there exists $[F]\in\Phi$ such that $\frak p_G=\frak p_F$. We observe by definition that $G,F$ are noetherian and so they are monoform modules by [LS, Proposition 4.7]. The fact that $\frak p_G=\frak p_F$ implies that $\overline{G}=\overline{F}$. Consequently $[G]=[F]$ by \cref{com}. Now, assume that $0\To N\To M\To K\To 0$ is an exact sequence of noetherian $A$-modules such that $\PSupp N,\PSupp K\subset \Phi$. Given $[H]\in\PSupp M$, we have $\frak p_H\in\Supp M$. Then $\frak p_H\in\Supp N$ or $\frak p_H\in\Supp K$ by \cref{f}. Usinig again \cref{f}, we have $[H]\in\PSupp N\cup \PSupp K$; and hence $[H]\in\Phi$. Conversely, suppose that $\frak p,\frak q$ are prime ideals of $A$ such that $\frak p\subset \frak q$ and $\frak p=\frak p_H$ for some $[H]\in\Phi$. Since $A/\frak q$ is a noetherian  $A/\frak p$-module, by Gruson theorem, $A/q$ has a finite filtration of submodules $$0=\subset L_1\subset\dots \subset L_{n-1}\subset A/\frak q$$  such that each $L_i/L_{i-1}$ is some quotient of finite direct sum of $H$. This implies that $A/\frak q\in\overline{C(H)}$ so that $[A/\frak q]\in\PSupp H\subset \Phi$. Therefore $\frak q\in\Phi$. 
\end{proof}
An immediate conclusion of the above proposition is a result due to Stanley and Wang [SW, Corollary 7.1]. 
\medskip
\begin{Corollary}\label{sertor}
Let $A$ be a noetherian ring. Then every nullity class of $\mod A$ is Serre.
\end{Corollary}
\begin{proof}
If $\cT$ is a nullity class of $\mod A$, then 
it follows from [LS, Theorem 6.8] that $\PSupp \cT$ is a closed and extension-closed subset of $\nPSpec A$. Then $\Supp\cT=\bigcup_{M\in\cT}\Supp M$ is a  specialization-closed subset of $\Spec A$ using \cref{f} and \cref{spec}. Thus $\cT$ is a Serre subcategory of $\mod A$ by [Ga].
\end{proof}

\medskip
\begin{Theorem}\label{homeom}
The function $p:\PSpec A\To\Spec A$; given by $[H]\mapsto \frak p_H$ is continuous. Moreover, if $A$ is notherian, then $p:\nPSpec A\To\Spec A$ is homeomorphism.  
\end{Theorem}
\begin{proof}
Using \cref{ann}, the function $p$ is well-defined. We show that $p^{-1}(\Phi)$ is a closed subset of $\PSpec A$ for every specialization closed subset $\Phi$ of $\Spec A$. Suppose that $[H]\in p^{-1}(\Phi)$ and so $\frak p_H=p([H])\in\Phi$. We prove that $\PSupp H\subset p^{-1}(\Phi)$. Given $[G]\in\PSupp H$, it follows from the definition that $G\in\overline{C(H)}$. Then $G$ contains a subobject $G_1$ which is a quotient of some premonoform  module $H_1$ isomorphic to $H$. Therefore using \cref{ann} and \cref{pp}, we have $\frak p_H=\Ann_AH_1\subset \Ann_A G_1=\Ann_A G=\frak q$. Since $\Phi$ is specialization closed, we deduce that $p([G])=\frak q\in\Phi$ so that $[G]\in p^{-1}(\Phi)$. To prove the second assertion, if $p([H])=p([G])$, then $\frak p_H=\frak p_G$. This implies that $\overline{H}=\overline{A/\frak p_H}=\overline{A/\frak p_G}=\overline{G}$. Now, \cref{com} forces $[H]=[G]$ so that $p$ is injective. For every $\frak p\in\Spec A$, we have $p([A/\frak p])=p$ so that $p$ is surjective.  
\end{proof}

\medskip
 \begin{Corollary}\label{extclos}
 Let $A$ be a noetherian ring and $0\To N\To M\To K\To 0$ be an exact sequence of noetherian $A$-modules. Then $\nPSupp M=\nPSupp N\cup \nPSupp K$. In particular, every closed subset of $\nPSpec A$ is extension-closed.  
 \end{Corollary}
 \begin{proof}
 It follows from \cref{f} and \cref{homeom} that $\nPSupp M=p^{-1}(\Supp M)=p^{-1}(\Supp N)\cup p^{-1}(\Supp K)=\nPSupp N\cup \nPSupp K$. The first assertion is clear by the first one.
 \end{proof}

 \begin{Corollary}\label{cooc}
 Let $A$ be a noetherian ring. Then the collection of closed  subsets of $\nPSpec A$ forms a topology.
 \end{Corollary}
 \begin{proof}
 Straightforward by \cref{extclos}.
 \end{proof}

 \begin{Corollary}
Let $A$ be a noetherian ring. Then the map $\cX\mapsto\nPSupp\cX$ establishes a one-to-one correspondence between localizing subcategories of $\cA$ and closed subsets of $\nPSpec A$. The inverse map is given $\Phi\mapsto \overrightarrow{\nPSupp^{-1}\Phi\cap\noeth A}$
\end{Corollary}
\begin{proof}
It follows from \cref{tors} that $\nPSupp^{-1}$ is a torsion class of $\Mod A$; and hence $\nPSupp^{-1}\Phi\cap\noeth A$ is a nullity class of $\mod A$. Thus \cref{sertor} implies that $\nPSupp^{-1}\Phi\cap\noeth A$ is a Serre subcategory of $\mod A$. Hence $\overrightarrow{\nPSupp^{-1}\Phi\cap\noeth A}$  is a localizing subcategory of finite type of $\Mod A$ by [Kr, Theorem 2.8]. Therefore $\overrightarrow{\nPSupp^{-1}\Phi\cap\noeth A}$  is a localizing subcategory of weakly finite type of $\cA$ by \cref{wft}. Now the result follows from \cref{classloc} and \cref{extclos}. 
 \end{proof}
 \medskip
Over a commutative noetherain ring $A$, every noetherian premonoform $A$-module is monoform, see [LS, Proposition 4.7]. We further have the following theorem.
 
\begin{Theorem}\label{spatom}
Let $A$ be a noetherian ring. Then there exists a bijective map  $\theta:\nPSpec A\to \ASpec A$ given by $[H]\mapsto \overline{H}$. Moreover, $\theta$ establishes a one-to-one correspondence between closed subsets of $\nPSpec A$ and open subsets of $\ASpec A$.    
\end{Theorem}
\begin{proof}
The map $\theta$ is well-defined by \cref{com}. To prove the second assertion,  it is clear that $\theta$ is surjective. Moreover, \cref{com} implies that $\theta$ is injective. If $\Psi$ is a closed subset of $\nPSpec\cA$ and $\Phi$ is an open subset of $\ASpec \cA$, then it is clear that $\Psi=\bigcup_{[H]\in\Psi}\nPSupp H$ and $\Phi=\bigcup_{\overline{H}\in\Phi}\ASupp H$.  Thus, it suffices to show that $\theta^{-1}(\ASupp M)=\nPSupp M$ and $\theta(\nPSupp M)=\ASupp M$ for any noetherian module $M$. Given $[H]\in\theta^{-1}(\ASupp M)$, we have $\overline{H}\in\ASupp M$. Hence there exists a monoform submodule $H_1$ of $H$ and a submodule $K$ of $M$ such that $H_1$ is a submodule of $M/K$. This implies that $H_1\in\sqrt{M}=\overline{C(M)}$ by \cref{sertor}. It now follows from \cref{com} that $H\in\overline{C(M)}$ and so $[H]\in\nPSupp M$. A similar proof shows that $\nPSupp M\subset\theta^{-1}(\ASupp M)$ and $\theta(\nPSupp M)=\ASupp M$.
\end{proof}

 \medskip
 
 A {\it perfect} $A$-complex $P_{\bullet}$ is a bounded complex $$P_{\bullet}:=0\To P_m\To P_{m-1}\To\dots\To P_n\To 0$$
formed by finitely generated projective $A$-modules $P_i$. We denote by $D(A)$, the derived category of $\Mod A$ and by $D_{\rm perf}(A)$,  the triangle subcategory of $D(A)$ consisting of $A$-complexes isomorphic to perfect $A$-complexes. 

 \begin{Definition}
 A full subcategory $\cX$ of a triangle category $\cT$ is said to be {\it thick} if $\cX$ is closed under taking direct summands and for every exact triangle $X\To Y\To Z\To \Sigma X$ in $\cT$, if two of the objects $X,Y,Z$ are in $\cX$, then so is the third.
  \end{Definition}
  
 Let $A$ be a noetherain ring. We denote by:

 $\bullet$ $\cL_{\rm thick}(D_{\rm perf}(A))$, the lattice of all thick subcategories of $D_{\rm perf}(A)$;
 
 $\bullet$ $\cL_{\rm null}(\mod A)$, the lattice of all nullity classes of $\mod A$;
 
 $\bullet$ $\cL_{\rm coh}(\mod A)$, the lattice of all coherent subcategories of $\mod A$;
 
 $\bullet$ $\cL_{\rm Serre}(\mod A)$, the lattice of all Serre subcategories of $\mod A$;
 
  $\bullet$ $\cL_{\rm narrow}(\mod A)$, the lattice of all narrow subcategories of $\mod A$;
 
 $\bullet$ $\cL_{\rm loc}(\Mod A)$, the lattice of all localizing subcategories of $\Mod A$;
 
  $\bullet$ $\cL_{\rm torsft}(\Mod A)$, the lattice of all torsion classes of finite type of $\Mod A$;
 
 $\bullet$ $\cL_{\rm torswft}(\Mod A)$, the lattice of all torsion classes of weakly finite type of $\Mod A$;
 
$\bullet$ $\cL_{\rm spcl}(\Spec A)$, the lattice of all specialization-closed subsets of $\Spec A$;

$\bullet$ $\cL_{\rm clos}(\nPSpec A)$, the lattice of all closed subsets of $\nPSpec A$;

$\bullet$ $\cL_{\rm cexc}(\nPSpec A)$, the lattice of all closed and extension-closed subsets of $\nPSpec A$;

$\bullet$ $\cL_{\rm sclos}(\nPSpec A)$, the lattice of all Serre closed subsets of $\nPSpec A$;

$\bullet$ $\cL_{\rm scsexc}(\nPSpec A)$, the lattice of all Serre closed and Serre extension-closed subsets of $\nPSpec A$;

$\bullet$ $\cL_{\rm open}(\ASpec A)$, the lattice of all open subsets of $\ASpec A$;

$\bullet$ $\cL_{\rm Zieg}(\Zg A)$, the lattice of all open subsets of $\Zg A.$

\medskip

  Hovey [Ho] established the following order-preserving maps between the lattices
  $$\xymatrix{
    \cL_{\rm spcl}(\Spec A) \ar@<1ex>[r]^{\Supp^{-1}} & \cL_{\rm loc}(\Mod A)\ar@<1ex>[r]^{\cap\mod A}\ar@<1ex>[l]^{\Supp} &  \cL_{\rm Serre}(\mod A) \ar@<1ex>[r]^{\cap\mod A}\ar@<1ex>[l]^{\langle\rangle_{\rm loc}} &   \cL_{\rm coh}(\mod A)\ar@<1ex>[r]^{f}\ar@<1ex>[l]^{\sqrt{}} &  \cL_{\rm thick}(D_{\rm perf}(A) \ar@<1ex>[l]^g} $$
in which  $f(\cC)=\{X_{\bullet}|\hspace{0.1cm} H_i(X_{\bullet})\in\cC$ for all $i\in\mathbb{Z}\}$ for every $\cC\in \cL_{\rm coh
}(\mod A)$ and  $g(\cT)$ is the coherent subcategory of $\Mod A$ generated by $\{H_i(X_{\bullet})|\hspace{0.1cm} X_{\bullet}\in\cT\}$ for every $\cT\in\cL_{\rm thick}(D_{\rm perf}(A).$   
\cref{spatom} induces following the order-preserving maps between lattices  
$$\xymatrix{
    \cL_{\rm clos}(\nPSpec A) \ar@<1ex>[r]^{\theta} & \cL_{\rm open}(\ASpec A)\ar@<1ex>[l]^{\theta^{-1}}} $$
 given by $\theta(\Psi)=\bigcup_{[H]\in\Phi}\ASupp H$ for every closed subset $\Psi\in \cL_{\rm closl}(\nPSpec A)$ and  $\theta^{-1}(\Phi)=\bigcup_{\overline{H}\in\Phi}\nPSupp H$ for every closed subset $\Phi\in\cL_{\rm open}(\ASpec A)$. \cref{spec} induces the following order-preserving maps between lattices  
$$\xymatrix{
  \cL_{\rm cexc}(\nPSpec A) \ar@<1ex>[r]^{p} &  \cL_{\rm spcl}(\Spec A)  \ar@<1ex>[l]^{p^{-1}}} $$
 given by $p(\Psi)=\{\frak p_H\in\Spec A|\hspace{0.1cm}  [H]\in\Phi\}$ for every closed subset $\Psi\in \cL_{\rm closl}(\nPSpec A)$ and  $p^{-1}(V)=\{[H]\in\nPSpec A|\hspace{0.1cm} \frak p_H\in V\}$ for every specialization-closed subset $V$ of $\Spec A$. 
 
 The Ziegler spectrum $\Zg \cA$ of a locally coherent Grothendieck category $\cA$ has been investigated in [H] which consists of all indecomposable injective objects in $\cA$.  For a finitely presented object $M$ in $\cA$, set $\cO(M)=\{E\in\Zg\cA|\hspace{0.1cm} \Hom(M,E)\neq 0\}$. A topology on $\Zg\cA$ formed by 
 \begin{center}
 $\{\cO(M)|\hspace{0.1cm} M$ is a finitely presented of $\cA\}$  
    \end{center}
as a basis of open subsets, is called the {\it Ziegler spectrum} of $\cA$. The Ziegler spectrum  $\Zg \Mod A$ is denoted by $\Zg A$. [K, Theorem 5.9] induces the following order-preserving maps between lattices  
$$\xymatrix{
  \cL_{\rm open}(\ASpec A) \ar@<1ex>[r]^{E} &  \cL_{\rm Zieg}(\Zg A)  \ar@<1ex>[l]^{E^{-1}}} $$
 given by $E(\Phi)=\{E(H)|\hspace{0.1cm} \overline{H}\in\Phi\}$ for every open subset $\Phi\in\Spec\cA$ and  $E^{-1}(\cO)=\{\overline{H}\in\ASpec A|\hspace{0.1cm} E(H)\in\cO\}$ for every open subset $\cO\in\cL_{\rm Zieg}(\Zg A)$.
 
 We now have the following corollary.

 \begin{Corollary}\label{lascor}
 Let $A$ be a noetherian ring. Then with the above notations, there are the following equalities and isomorphisms of lattices:
\[ 
\begin{array}{l}
 \cL_{\rm Zieg}(\Zg A)\stackrel{E^{-1}}\cong\cL_{\rm open}(\ASpec A)
 \overset{\theta^{-1}}\cong \cL_{\rm clos}(\nPSpec A)\\
  \overset{(1)}=\cL_{\rm cexc}(\nPSpec A)\overset{(2)}=\cL_{\rm sclos}(\Spec A)
  \overset{(3)}=\cL_{\rm cexc}(\nPSpec A)\\\overset{p}\cong \cL_{\rm spcl}(\Spec A)\overset{Supp^{-1}}\cong \cL_{\rm loc}(\Mod A)\overset{(4)}=\cL_{\rm torsft}(\Mod A)
  \\\overset{(5)}=\cL_{\rm torswft}(\Mod A)\overset{\cap\mod A}\cong \cL_{\rm Serre}(\mod A)
  \overset{(6)}=\cL_{\rm null}(\mod A)\\\overset{(7)}=\cL_{\rm narrow}(\mod A)\overset{(8)}=\cL_{\rm coh}(\mod A)\overset{f}\cong \cL_{\rm thick}(D_{\rm perf}(A)).
 \end{array}
 \]
  \end{Corollary}
\begin{proof}
The map $E^{-1}$ is isomorphism by [K, Theorem 5.9]. The map  $\theta^{-1}$ is isomorphism by \cref{spatom}. The equality (1) holds by \cref{extclos}. The equalities  (2) and (3) holds by \cref{comclex}. The map $p$ is isomorphism by \cref{spec} and \cref{homeom}. The map $\Supp^{-1}$ defined as $\Supp^{-1}(V)=\{M\in \Mod A|\hspace{0.1cm} \Supp M\subset V\}$ for every $V\in\cL_{\rm loc}(\Mod A)$, is isomorphism by [T, Theorem 3.6]. To show the equality (4), since every $A$-module in a localizing subcategory $\cX$ of $\Mod A$ is the direct union of its noetherain subobjects, we conclude that $\cX=\overrightarrow{\cX\cap\mod A}$. Hence, it follows from [Kr, Theorem 2.8] that $\cX$ is of finite type of $\Mod A$ and so $\cL_{\rm loc}(\Mod A)\subset \cL_{\rm torsft}(\Mod A)$. Conversely, if $\cY\in  \cL_{\rm torsft}(\Mod A)$, we have $\cY=\overrightarrow{\cY\cap\mod A}$ by \cref{wft}. We observe that $\cY\cap\mod A$ is a nullity class of $\mod A$ and so it is Serre subcategory of $\mod A$ by \cref{sertor}. Now, using again [Kr, Theorem 2.8], the subcategory $\cY$ of $\Mod A$ is localizing so that $   \cL_{\rm torsft}(\Mod A)\subset\cL_{\rm loc}(\Mod A)$. The equality (5) is straightforward by \cref{wft} and [Kr, Theorem 2.8]. The map $\cap \mod A$ is isomorphism by [T, Theorem 3.6].  The equality (6), (7) and (8) hold by [SW, Corollary 7.1]. The map $f$ is isomorphism by [T, Theorem 3.6]. 
\end{proof}


\end{document}